\documentclass[12pt]{article}
\usepackage{fullpage}
 \usepackage{amsmath}
\usepackage{amssymb}
\usepackage{amsbsy}
\usepackage{amsfonts}
\def\rit{{\mathbb R}}

\def\zit{{\mathbb Z}}
\def\tit{{\mathbb T}}

\def\eps{\varepsilon}

\def\beq{\begin{equation}}
\def\eeq{\end{equation}}

\def\Re{{\rm Re}  \,}
\def\Im{{\rm Im}  \,}
\newtheorem{theo}{Theorem}
\newtheorem{prop}[theo]{Proposition}
\newtheorem{lem}[theo]{Lemma}
\newtheorem{rem}[theo]{Remark}

\title
{Transverse nonlinear instability for two-dimensional dispersive models }
\author{F. Rousset\footnote{University of Nice-Sophia Antipolis, Laboratoire J.A. Dieudonn\'e, 
06108 Nice cedex 2 },\,\, 
N. Tzvetkov\footnote{University of Lille 1, Laboratoire P. Painlev\'e, 59655 Villeneuve d'Ascq cedex}}
\date{}
\begin{document}
\maketitle
\begin{abstract}
We present a method to prove nonlinear instability of solitary waves  in dispersive models. 
Two examples are analyzed: 
we prove the nonlinear long time instability of the KdV solitary wave (with respect 
to periodic transverse perturbations) under a KP-I flow
and the transverse nonlinear instability  of solitary waves for the cubic nonlinear Schr\"odinger equation.
\end{abstract}
\section{Introduction}
There are many results (both theoretical and numerical) dealing with detecting unstable modes 
of dispersive equations linearized around soliton like structures. 
However, in most of these cases it is not clear whether one has indeed a nonlinear instability for a flow
of the full nonlinear problem. The main reason is the lack of understanding of the whole spectrum
of the linearized problem. The goal of this paper is to present a method showing how only a partial information
about the spectrum of the linearized operator together with a suitable nonlinear analysis may indeed 
give the proof of the nonlinear instability in the presence of an unstable mode.
Our first example is the nonlinear long time instability of the KdV solitary wave (with respect 
to periodic transverse perturbations) under a KP-I flow.
We also prove a nonlinear instability result for the cubic nonlinear Schr\"odinger equation.
We believe that the method presented here could be useful in the contexts of other dispersive equations.
\\

Consider the Kortweg- de Vries (KdV) equation
\beq\label{KdV}
u_{t} +u u_{x} + u_{xxx} =0,
\eeq
$u:\mathbb{R}^2\rightarrow \mathbb{R}$,
which is an asymptotic model, derived from the free surface Euler equation, 
for the propagation of long one-directional small amplitude surface waves.
A famous solution of (\ref{KdV}) is the solitary wave solution, given  by 
$$
u(t,x) = Q(x-t), \quad Q(x) = 3\,  \mbox{sech}^2\Bigl({ \frac{x}{2}}\Bigr).
$$
Observe that $u(t,x)$ corresponds to the displacement of the profile $Q$ from left to the right with speed
one. One also has the solution
\begin{equation}\label{speed}
u_{c}(t,x) =  c \,Q( \sqrt{c}(x- c t )), \quad c>0
\end{equation}
which corresponds to a solitary wave with a positive speed $c$.
 
A very natural question concerning the relevance of the solution $Q(x-t)$ is its stability with respect
to small perturbations. It is evident that the usual Lyapounov stability cannot hold because of the 
translation invariance of the problem. More precisely for $c$ close to one $c \,Q( \sqrt{c}\,x)$ is close to
$Q(x)$ while for $t\gg 1$ ($t\sim|c-1|^{-1}$) the corresponding solutions of the KdV equation $u(t,x)$ and
$u_{c}(t,x)$ separate from each other at distance independent of the smallness of $c-1$.
However, the solution $u_{c}(t,x)$ remains close to the spatial translates of $Q$ and thus orbital stability
of $Q$ under the flow of KdV is not excluded. It is known since the seminal paper of Benjamin~\cite{Ben}
that $Q$ is orbitally stable in the energy space $H^1(\mathbb{R})$ (we call $H^1(\mathbb{R})$ the energy space
since this is the natural space induced by the Hamiltonian structure of (\ref{KdV})). Here is the precise 
statement.
\begin{theo}\label{th1}
For every $\varepsilon>0$ there exists $\delta>0$ such that if the initial data 
$$
u|_{t=0}=u_{0}\in H^1(\mathbb{R})
$$ 
of the KdV equation (\ref{KdV}) satisfies
$$
\|u_0-Q\|_{H^1(\mathbb{R})}<\delta
$$
then the corresponding solution $u$ (which is well defined thanks to \cite{KPV}) satisfies
$$
\sup_{t\in \mathbb{R}}\inf_{a\in \mathbb{R}}\|u(t,\cdot)-Q(\cdot-a)\|_{H^1(\mathbb{R})}<\varepsilon\,.
$$
\end{theo}
Let us notice that the phase space $ H^1(\mathbb{R})$ may be replaced by $ L^2(\mathbb{R})$ (see \cite{MV}).

In \cite{KP}, Kadomtsev-Petviashvili studied weak transverse perturbation of the KdV flow and derived 
the following two-dimensional models
\begin{equation}\label{kp}
u_{t} +u u_{x} + u_{xxx}  \pm\partial_{x}^{-1} u_{yy}  = 0\,.
\end{equation}
Equation (\ref{kp}) with sign $+$ is called the KP-II equation while (\ref{kp})
with sign $-$ is the KP-I equation. Let us observe that in the derivation of the model, the signs
vary in front of the $u_{xxx}$ term  and correspond  to different surface tensions.
However from mathematical view point the study of the models with signs varying in front of 
$u_{xxx}$ is equivalent to the study of the models
with signs varying in front of $\partial_{x}^{-1} u_{yy}$ by the variable change
$
u(t,x,y)\mapsto -u(-t,x,y).
$
The anti-derivative $\partial_{x}^{-1}$ is defined on functions which have, in a suitable sens, a zero $x$
mean value. 

Let us observe that $Q(x-t)$ is a solution of both equations (\ref{kp}). It is conjectured in \cite{KP}
that $Q(x-t)$ is stable under the KP-II flow and unstable under the KP-I flow. Of course this conjecture 
is very vague since one should precise the stability notion and the spatial domain for $x$, $y$.
In  \cite{Pego}, all possible unstable modes of the linearized
 equation are described and in particular it is shown that the linearization
about $Q$ of the KP-I flow is unstable and the linearization of the  KP-II flow is spectrally stable.
In this paper, we show that the spectral  instability result of \cite{Pego} indeed implies the
nonlinear instability in the case of the KP-I equation for solutions periodic
in the $y$ variable. This  result is actually not new since the equation being
completely integrable (having a Lax pair representation), 
the instability can be shown by exhibiting explicit solutions (see Zakharov~\cite{Z}). 
Nevertheless, we believe that our method inspired from the work of 
Grenier~\cite{Grenier} in fluid mechanics to prove that spectral instability implies nonlinear 
instability which does not use the complete integrability  is interesting and can be applied 
to many other dispersive equations. As an illustration, we shall also study below a transverse 
instability of the two-dimensional cubic nonlinear Schr\"odinger equation which is not completely integrable.
    
The global well-posedness of the KP-I equation in the setting ${\mathbb R}\times{\mathbb T}$ 
was recently obtained by Ionescu-Kenig \cite{Kenig} in a space
which moreover  contains the solitary wave $Q$ and hence,  we state our result in the 
context of Ionescu-Kenig's theorem. In general it is difficult to get  nonlinear instability 
results in natural  energy  norms like  $L^2$ or $H^1$   for conservative equations
due to the presence of strong  nonlinearities.  Here we shall use the general setting developed by Grenier in 
\cite{Grenier} in the context of the Euler equation which relies on the possibility of constructing an 
high order  approximate solution more accurate that the  only linear approximation. For other methods,  
we refer to \cite{FSV}, \cite{GS}.

We consider thus the KP-I equation
\beq\label{KP}
u_{t} +u u_{x} + u_{xxx}  -  \partial_{x}^{-1} u_{yy}  = 0
\eeq
for $(x, y ) \in \rit\times \tit_{L}$ where $\mathbb{T}_{L}$ is the flat torus 
$\mathbb{R}/ 2\pi L \mathbb{Z}$. 
As mentioned above, a special solution of this equation is given by the KdV soliton
$Q(x-t)$.
Since we are interested in the stability of the soliton for \eqref{KP}, it is more convenient to go 
into a moving frame  i.e. to change $x$ into $x-t$ and to study the equation
\beq\label{KPm}
u_{t} - u_{x} +u u_{x} + u_{xxx}  -  \partial_{x}^{-1} u_{yy}  = 0, \quad (x, y) \in \rit \times\tit_{L}
\eeq
so that $Q(x)$ is now a stationary solution of \eqref{KPm}. 
Note that we can  always 
change space and time scales to reduce the study  of the stability of $u_{c}$, given by (\ref{speed}) 
to the study of the stability of $Q$ for \eqref{KPm}. 
Nevertheless, since we are in a bounded domain in $y$, the scaling 
changes the size of the domain, this is why we keep the parameter $L$ in our study. 

As established in \cite{Kenig}, the Cauchy problem for \eqref{KP} or equivalently \eqref{KPm} is  
globally well-posed for data in  the space $Z^2( \rit \times \tit_{L})$ defined by 
$$ 
Z^2(\rit \times \tit_{L}) = 
\Bigl\{ u, \, \| \hat{u} (\xi, k) \, ( 1 +|\xi|^2 + |k/\xi |^2 ) \|_{ L^2( \rit \times \zit)} 
<+ \infty \Bigr\},
$$
where $\hat{u}(\xi, k)$ is the Fourier transform of $u$ : 
$$ 
\hat{u}(\xi, k)=  \frac{1}{2\pi L} \int_{-\infty}^{\infty}
\int_{0}^{2\pi L} e^{ -i \big(  x \, \xi + \frac{ y\, k} {L} \big)} u(x,y) \, dydx.
$$
If  $u\in Z^2$, this means that $ u,\, u_{x}, \, u_{xx}$ and 
$\partial_{x}^{-1}u_{y}, \, \partial_{x}^{-2}u_{yy}$ are in $L^2$,
where $\partial_{x}^{-1}$ is defined in the natural way via the Fourier
transform for functions $u\in L^2$
such that $\xi^{-1}\hat{u}(\xi, k)\in L^2$. Moreover, the propagation
of $H^s$ regularity holds: if  $u_{0} \in H^s \cap Z^2 $ for $s>7$, then
$ u(t) \in  H^s \cap Z^2 $ for every $t>0$.  Note that since $Q$ does not depend of $y$, we have $Q \in Z^2$.
The first goal of this paper is to prove the following orbital instability result.
\begin{theo}\label{maintheo} Assuming that $L>4/\sqrt{3}$, then for every $s\geq 2$, there exists 
$\eta>0$ such that  for every $\delta>0$,   there exists $u_{0}^\delta\in Z^2\cap H^s$  and a time 
$T^\delta \sim |\log \delta |$ such that
$$ 
\|u_{0}^\delta - Q \|_{H^s(\rit \times \tit_{L})} < \delta
$$
and the solution  $u^\delta $   of \eqref{KPm} with initial value $u_{0}^\delta$ satisfies 
$$ 
\inf_{a\in \rit}\|u^\delta(T^\delta,\cdot)- Q(\cdot-a)\|_{L^2(\rit \times \tit_{L})} \geq \eta.
$$ 
\end{theo}
\begin{rem}
If $u(t,x,y)$ is a solution of the KP equation (\ref{kp}) then so is $u_{\lambda}$ defined by
$$
u_{\lambda}(t,x,y)=\lambda^2 u(\lambda^3 t,\lambda x,\lambda^2 y).
$$
Thus in the context of (\ref{KP}) 
solutions of period one in $y$ transform into solutions of period $\lambda^{-2}$ and solitary
waves of speed $c$ transform into solitary waves of speed $\lambda^2 c$.
Consequently Theorem~\ref{maintheo} implies that 
if we fix period one perturbations in $y$ then one 
needs to consider solitary waves of sufficiently large speed
to get the instability.
Let us also remark that the restriction $L>4/\sqrt{3}$ in Theorem~\ref{maintheo} is imposed by the 
spectral considerations of \cite{Pego} and is needed for the existence of unstable eigenmodes.
It would be interesting to decide what happens for $L\leq4/\sqrt{3}$ (or equivalently for small
speed solitary waves for period one perturbations).
\end{rem}
\begin{rem}
Let us recall that a three dimensional analogue of (\ref{KP}) 
\begin{equation}\label{3D}
u_{t} +u u_{x} + u_{xxx}  -  \partial_{x}^{-1} (u_{yy}+u_{zz})  = 0, \quad (x,y,z)\in \rit^3
\end{equation}
has solutions blowing up in finite time (see \cite{Liu} and also \cite{saut}) and thus for the 
three dimensional versions
of the KP-I equation a stronger form of the instability appears. It is however an open problem to prove the
existence of blow-up solutions for (\ref{3D}) with $u$ periodic in $y$, $z$.
\end{rem}
Let us outline the main steps of the proof  Theorem~\ref{maintheo}.
First, we need to use the result of \cite{Pego} concerning the existence of unstable
eigenmodes for the linearized about $Q$ operator. 
Next,  following the idea of
Grenier \cite{Grenier}, we perform the construction of an
approximate solution. 
The approximate solution is defined iteratively. At the first step we put the unstable
eigenmode.
At each further step, we get linear problems with source terms involving the previous iterates
(the procedure is closely related to the Picard iteration).
We need to control precisely  the eventual growth in time of each iterate. By applying a Laplace transform, 
we reduce the matters to showing estimates on a resolvent equation which
are uniform on some  straight line $\lambda = \gamma + i \tau, \, \tau \in \mathbb{R}.$ 
For  bounded  frequencies  (i.e. $|\tau|$ bounded), a classical ODE argument combined with the absence 
of unstable modes coming from \cite{Pego} 
suffices to get the needed bound.  The main difficulty is  to get uniform   resolvent estimates for 
large $\tau$. They will 
result from conservation (or almost conservation) laws.
We finally perform an energy estimate to the nonlinear problem to show that the constructed approximate
solution is indeed close to the actual solution for suitable time scales. This in turn implies 
the nonlinear instability claimed in Theorem~\ref{maintheo}.
\\

The second example that we  consider in this paper is  the 
two-dimensional  Nonlinear Schr\"odinger equation (NLS) 
\beq
\label{NLS0}
 i v_{t} + \Delta_{x,y} v   + |v|^2 v = 0,
\eeq
where $v$ is a complex valued function.
A famous solution of this equation is the solitary wave
 $ Q(x)  e^{it}$ with $Q$ given by 
$$
Q(x)=\frac{\sqrt{2}}{\mbox{ch} (x)}\,.
$$ 
This solitary wave is  orbitally stable when submitted to one-dimensional perturbations i.e. 
perturbations which depend on $x$ only (see \cite{CL}). Here orbital stability means that 
$$
\forall\, \varepsilon>0,\,\, \exists\,\, \delta>0\,:\,
\|v(0,\cdot) - Q \|_{H^1(\mathbb{R})}<\delta
\implies
\inf_{a\in \mathbb{R}, \, \gamma \in \mathbb{R} }  
\|v(t,\cdot) - e^{i\gamma} e^{it} Q(\cdot-a) \|_{H^1(\mathbb{R})}<\varepsilon\, .
$$
We shall prove that, similarly to the KdV soliton as a solution of the KP-I equation, 
this stationary solution of (\ref{NLS}) which is orbitally stable when submitted 
to one-dimensional perturbation is nonlinearly unstable when it is submitted to two-dimensional
perturbation.  As previously, it is more convenient to set $ v = e^{it} u$ and to study
 the equation 
\beq\label{NLS}
i u_{t} + \Delta u - u + |u|^2 u = 0,
\eeq
for $(x,y)\in \mathbb{R}\times \mathbb{T}_{L}$. A stationary solution of this equation
is now given by the ground state  $Q(x)$. Since the solitary waves modelled on 
$Q(x)$  for \eqref{NLS0} are given by 
$$
u_{\lambda}(t,x) = \lambda Q(\lambda x) e^{i \lambda^2 t}
$$
we can always reduce by scaling  the study of the stability  of  $u_{\lambda}$ to the study of 
the stability of $Q$ in \eqref{NLS}, but it is again important to keep $L$ as a parameter.
Here is our result. 
\begin{theo}\label{theoschro} 
There exists $L_{0}$ such that for $L\geq L_{0}$,   for every $s\geq 2$, there exists 
$\eta>0$ such that  for every $\delta>0$,   there exists $u_{0}^\delta\in  H^s$  and a time 
$T^\delta \sim |\log \delta |$ such that
$$ 
\|u_{0}^\delta - Q \|_{H^s(\rit \times \tit_{L})} < \delta
$$
and the solution  $u^\delta $   of \eqref{NLS} with initial value $u_{0}^\delta$  
belongs to $\mathcal{C}([0, T^\delta], H^s)$ and satisfies 
$$ 
\inf_{a\in \rit, \, \gamma \in \mathbb{R} }\|u^\delta(T^\delta,\cdot)-
Q(\cdot-a)e^{i\gamma}\|_{L^2(\rit \times \tit_{L})} 
\geq \eta.
$$ 
\end{theo}
\begin{rem}
Let us observe that the cubic two-dimensional NLS
is not known to be integrable (in the sense of Lax pairs representation) 
and thus it is hard to expect that the instability result presented in 
Theorem~\ref{theoschro} can be displayed by an explicit family of solutions in the spirit of \cite{Z}.
\end{rem}
\begin{rem}
It is likely that the method presented here may be applied to the following
two dimensional perturbation of the Boussinesq equation
\begin{equation}\label{bous}
u_{tt}+(u_{xx}+u^2-u)_{xx}\pm u_{yy}=0\,.
\end{equation}
The stability, for suitable values of the propagation speed, with respect to
one-dimensional perturbations of the soliton of the Boussinesq equation is
obtained in \cite{BS}. The analysis for an unstable mode in 2D in the context
of (\ref{bous}) is essentially the same as the corresponding analysis for the
KdV soliton as a solution of KP (see \cite{BBD}). One thus may perform the
analysis of \cite{Pego} (see also the Appendix of this paper) combined with
the nonlinear analysis of this paper to get statements in the spirit of
Theorems~\ref{maintheo}, \ref{theoschro} for equation (\ref{bous}).
\end{rem}
The assumption $L\geq L_{0}$ in Theorem~\ref{theoschro} is used to get the spectral instability of the solitary wave.  
A difference with Theorem~\ref{maintheo} is that for the two-dimensional Schr\"odinger equation 
in $\mathbb{R} \times \mathbb{T}$ a global existence result of large data strong solutions is 
not known so that
Theorem~\ref{theoschro} contains the fact that our unstable solution $u^\delta$  remains well-defined  
on a sufficiently long time scale.  
In fact, small data global existence for (\ref{NLS}), posed on $\mathbb{R} \times \mathbb{T}$, 
is obtained in \cite{TT}. For general large data we may not have the global existence for 
(\ref{NLS}), posed on $\mathbb{R} \times \mathbb{T}$, since one may localize the well-known 
explicit blow-up solution for the cubic NLS on ${\mathbb R}^2$ (see \cite{BGT} for details on this argument).

The rest of the paper is organized as follows. In the next section,  we give a detailed proof of 
all the steps of the proof of Theorem~\ref{maintheo}. Then 
we give a less detailed  proof of Theorem~\ref{theoschro} since the method is the same. 
Finally, the  appendix is devoted to the linear instability results. 
\\

{\bf Acknowledgement.} We are grateful to Robert Pego for pointing out to us the reference~\cite{Z}.
A previous version of this text, before we were aware of the Zakharov work \cite{Z}, was posted
to the arxiv of preprints on December 2006.
\section{Proof of Theorem \ref{maintheo}}
\subsection{Existence of a most  unstable eigenmode}
The linearized equation about the soliton  $Q$ reads
\beq\label{eqlin}
u_{t} +A u = 0, \quad Au =-u_{x}+ (Qu)_{x} + u_{xxx}  - \partial_{x}^{-1}\partial_{yy}u,\quad 
(x,y)\in\mathbb{R}\times \mathbb{T}_{L}\,.
\eeq
This last linear equation can be solved, for instance by a classical energy method, for initial data
in $H^s$ such that its anti-derivative exists.     
The main result of \cite{Pego} is the characterization of all the unstable eigenmodes associated to $A$. 
An unstable eigenmode is a solution of (\ref{eqlin}) under the form
$$  
\varphi_{\sigma, k }(t,x,y)= e^{\sigma t } e^{\frac{i k y} {L} } V(x), 
$$
with $\Re \sigma >0 $, $ V \in L^2(\mathbb{R})$. The result of \cite{Pego} adapted to  our framework reads : 
\begin{theo}[\cite{Pego}]\label{Pego}
There exists unstable eigenmodes   if and only if $L > \frac{4}{\sqrt{3}}$.  Moreover,   for an 
unstable eigenmode, $\sigma$ and $k\in \mathbb{Z}$  are parametrized by
\beq\label{instab}
2 \sigma = \mu(\mu-1)(2- \mu), \quad k = \frac{\sqrt{3}  \,L}{ 4 } \,\mu(2-\mu),\quad \mu \in (1,2)
\eeq
and there exists $g \in H^\infty(\mathbb{R})$ such that
\beq\label{propV}
V= g_{xx}.
\eeq
\end{theo}
For the sake of completeness, we recall the main steps of the proof of this result  in the Appendix.
         
Note that for $\mu \in (0,2)$, $\mu(2-\mu) \in(0,1)$ hence one can find an integer such that 
$k = \frac{\sqrt{3}  \,L }{ 4 } \,\mu(2-\mu)$ only if $\frac{\sqrt{3}  \,L}{ 4 }>1$. Moreover, for $L$ fixed, 
there is only a finite number of $k$ which verify this property, this allows us to choose $\sigma_{0}$ and 
$k_{0}$ such that $\varphi_{\sigma_{0}, k_{0}}$ is the most unstable eigenmode i.e.
$$ 
\sigma_{0} = \sup \Bigl\{ \sigma, \, (\sigma , k) \mbox{ verifying \eqref{instab} }\Bigr\}
$$
and $k_0$ is the corresponding integer such that \eqref{instab} holds with $(\sigma,k)=(\sigma_0,k_0)$.
Let us define  
$$
u^{0}(t,x,y)\equiv \varphi_{\sigma_{0}, k_{0}}(t,x,y)+ \overline{ \varphi_{\sigma_{0}, k_{0}}}(t,x,y)\,.
$$
To prove Theorem \ref{maintheo}, we shall use $Q+\delta u^0$ as an initial data for \eqref{KPm}.
As remarked before, we have  $Q\in Z^2 \cap H^s$ for every $s$, but  thanks to \eqref{propV} in 
Theorem~\ref{Pego}, we also have that $ u^0\in Z^2 \cap H^s$ consequently, thanks to the result of 
\cite{Kenig} there is a unique global  solution  $u^\delta $ of $\eqref{KPm}$  in $ Z^2 \cap H^s$ 
with initial value $Q+  \delta u^0$.  So the only problem that remains is to estimate from below
$$ 
\inf_{a\in \rit}\|u^\delta(T^\delta,\cdot)- Q(\cdot-a)\|_{L^2(\rit \times \tit_{L})}.
$$ 
Towards this, we shall use the method of \cite{Grenier} which relies on the construction 
of an high order unstable  solution. This is the aim of the next section.      
\subsection{Construction of an high order unstable approximate solution}
Let us set $ v = u^\delta - Q$, then $v$ solves
\beq\label{nonlin}
v_{t} + A v = - v v_{x}\,.
\eeq
We define $V_{K}^s$ as the space :
$$ 
V_{K}^s =\Bigl\{ u, \, u= \sum_{j \in \frac{k_{0}}{ L} \mathbb{Z}, \, |j \, L/k_{0}| \leq K} 
u_{j}(x) e^{  i  j  y  },\quad \, u_{j} \in H^s (\rit) \Bigr\} 
$$
and we define a norm on $V_{K}^s$ by
$$ 
|u|_{V_{K}^s} = \sup_{j} |u_{j}|_{s}
$$ 
where $|\cdot |_{s}$ is the standard $H^s(\rit)$ norm. 
Let us notice that  $u^{0} $ is such that $u^{0} \in V^{s}_{1}$ for all $s\in\mathbb{N}$.
Following the strategy of \cite{Grenier}, for $s\gg 1$, we look for an high order solution under the form :
\begin{equation}\label{u_app}
u^{ap}= \delta\Bigl( u^0  + \sum_{ k=1}^M \delta^k u^k \Bigr), \quad u^k \in  V^{s-k }_{ k +1}
\end{equation}
such that $u^k_{/t=0}=0$ and $M\geq 1$ is to be fixed later.
Once the value of $M$ is fixed, then we fix the integer $s$ so that $s>M$. 
By plugging the expansion in \eqref{nonlin}, 
cancelling the terms involving $\delta^k$, $1\leq k\leq M+1$,
we choose $u^k$ so that $u^k$ solves the problem
\beq\label{uk}
\partial_{t} u^k + A u^k = -\frac{ 1}{ 2 } 
\Bigl( \sum_{j+l= k-1} u^j u^{l} \Bigr)_{x}, \, \quad u^k_{/t=0} = 0.
\eeq  
The main point in the analysis of $u^{ap}$ is the following estimate.
\begin{prop}\label{propuk}
Let $u^k$ the solution of \eqref{uk},  if $s-k\geq 1$,  we have the estimate:
\beq\label{estuk} 
|u^k(t)|_{V^{s- k }_{k+1}} \leq C_{k,s}e^{ (k+1 ) \sigma_{0}t}, \, \forall t \geq 0.
\eeq 
\end{prop}
The proof of the proposition will follow easily by induction from the following theorem.
\begin{theo}\label{theolinsource}
Consider  the solution $u$ of the linear problem
\beq\label{linsource}
\partial_{t} u + A u = F_{x}, \quad u_{/t=0}= 0
\eeq
with a source term $F \in V_{K}^{s+ 1} $   with
\beq\label{F}
|F(t)|_{V_{K}^{s+ 1}} \leq C^F_{K,s} e^{ \gamma t}, \quad \gamma \geq 2\sigma_0
\eeq 
then $u$ belongs to $ V_{K}^{s}$ and satisfies the estimate
\beq\label{ulin}
|u(t)|_{V_{K}^{s}} \leq C_{K, s} e^{ \gamma t}, \quad \forall t \geq 0.
\eeq 
\end{theo}
We first observe that under our hypothesis on $F$ the solution of (\ref{theolinsource})
is well-defined and the only point is to  prove the quantitative bound (\ref{ulin}). 
The estimate  \eqref{ulin}, relies on  the fact that on $V_{K}^s$, 
the real part of the  spectrum of the operator $-A$ is bounded by $\sigma_0$. 
Nevertheless for such  a dispersive operator,  there is no general theory  to convert an information 
on the position of the spectrum into an estimate on the semi-group   like it is the case for example 
for sectorial operators. To get the result,  we  need to estimate the resolvent of $-A$ on $V_{K}^s$.          
At first, we can perform some reductions on the problem. Indeed, since $F$ has a finite number of 
Fourier modes, we can expand $u$ in Fourier modes and hence we only need to study the problem
\beq\label{ulink}
\partial_{t}v + A_{j}v =  \partial_{x}F_{j}(t,x),  \, \quad v_{/t=0}= 0 
\eeq
where 
\begin{equation}\label{Aj}
A_{j} v =  - v_{x} + ( Q v)_{x} + v_{xxx} + j^2 \partial_{x}^{-1} v,
\end{equation}
$j \in \frac{k_{0}}{ L} \mathbb{Z}$,
$|j \, L/k_{0}| \leq K$ and $v(t,x)=u_{j}(t,x)$, and to establish that $v$
satisfies 
$$
|v(t)|_{s}\leq C_{j,s}e^{\gamma t}
$$
under the assumption
\begin{equation}\label{Fbis}
|F_j(t)|_{s+1}\leq C_{j,s}e^{\gamma t}\,.
\end{equation}
In what follows, we fix $\gamma_{0}$ such that  $\sigma_0 <\gamma_{0}<\gamma$ and we shall use the 
Laplace transform. For $T>0$, we first   introduce $G$ such that
$$ 
G=0, \, t<0, \quad G=0, \, t>T, \quad G= F_j,  \, t \in [0, T]
$$ 
and we notice that the solution of 
$$
\partial_{t} \tilde{v} + A_{j } \tilde{v} = G_{x}, \quad \tilde{v}_{/t=0}=0
$$
coincides with $v$ on $[0,T]$ so that it is sufficient to study $\tilde{v}$. Next, we set
$$ 
w(\tau,x)= \mathcal{L} \tilde{v}(\gamma_{0}+ i \tau), \quad 
H(\tau, x) = \mathcal{L}G ( \gamma_{0} + i \tau ), \quad (\tau, x) \in \mathbb{R}^2
$$
where $\mathcal{L}$ stands for the Laplace transform in time :
$$ 
\mathcal{L}f(\gamma_{0}+ i  \tau) = \int_{\mathbb{R}} e^{- \gamma_{0} t - i \tau\, t} 
f(t) {\bf 1}_{t \geq 0}\, dt.
$$
We get that $w$ solves the resolvent equation
\beq\label{w}
(\gamma_{0} + i \tau ) w + A_{j} w = H_{x}.
\eeq
In the sequel, for complex  valued functions depending on $x$, we define
$$ 
(f, g ) \equiv  \int_{\mathbb{R}} f(x) \overline{g}(x)\, dx, \quad 
|f|^2\equiv \|f\|_{L^2(\mathbb{R})}^2 = (f,f), \quad |f|_{s}^2\equiv \|f\|_{H^s(\mathbb{R})}^2 
=\sum_{0\leq m  \leq s} |\partial_{x}^m f|^2 .
$$
Towards the proof of Theorem \ref{theolinsource}, we first need to study \eqref{w}. 
Our main estimate on the resolvent will be
\begin{theo}[Resolvent Estimates]\label{resolvant}
Let $s\geq 1$ be an integer.
Let $w(\tau)$ be the solution of \eqref{w} for $j$, $|j | \leq k_{0} K/L $, then there exists 
$C(s,\gamma_{0}, K)>0$ such that  for every $\tau$, we have the estimate
\beq\label{estresolvant}
|w(\tau)|^2_{s} \leq C(s,  \gamma_{0}, K) |H(\tau)|^2_{s+1}.
\eeq
\end{theo}
\subsubsection{Proof of Theorem \ref{resolvant}}
 We shall split the proof  in various lemmas. To estimate $w$, we shall deal differently with  large 
and bounded frequencies.           
\begin{lem}\label{HF}
There exists $M>0$ (which depends on $K$) and $C(s, \gamma_{0}, K)$
  such that for $|\tau|\geq M$, we have the estimate
\beq\label{estHF}
|w(\tau)|_{s}^2 \leq C(s, \gamma_{0}, K) |H(\tau)|_{s+1}^2.
\eeq
\end{lem}
\subsubsection{Proof of Lemma \ref{HF}}
We first prove \eqref{estHF} for $s=1$.       
Note that the equation \eqref{w} can be rewritten as
\beq\label{wL}
(\gamma_{0} + i \tau ) w - ( \mathcal{L} w )_{x} + j^2 \partial_{x}^{-1} w = H_{x}
\eeq
where $\mathcal{L}$ is defined by
$$ 
\mathcal{L}w =  - w_{xx} - Q w + w.
$$
Note that $\mathcal{L}$ is a self-adjoint operator in $L^2$ which is very useful in the proof of the stability of the 
soliton for the KdV equation.   Since it is self-adjoint, the spectrum is real. Moreover, since  
$Q$  goes to zero exponentially fast, the essential spectrum of $\mathcal{L}$ is in  $ [1, + \infty)$.  
For $\lambda< 1$ there are only eigenvalues of finite multiplicity.  Finally 
by Sturm-Liouville theory, since $Q_{x}$ is in the kernel of $\mathcal{L}$ and has only one zero,   
we get that $\mathcal{L}$ has  only one negative eigenvalue. Moreover,  $0$ is a simple eigenvalue. 
Consequently we can define an orthogonal decomposition: 
\beq\label{dec}
w = \alpha(\tau) \varphi_{-1} + \beta(\tau) \varphi_{0} + w_{\perp}
\eeq
where 
\beq\label{decprop}
\mathcal{L} \varphi_{-1} = \mu \varphi_{-1}, \, \mu <0, \,  \mathcal{L} \varphi_{0} = 0, \, 
( \mathcal{L}w_{\perp}, w_{\perp}) \geq c_{0} |w_{\perp}|^2, \, c_{0}>0.
\eeq
Note that the eigenvectors $\varphi_{-1}$ and $\varphi_{0}$ are smooth.
The important role of $\mathcal{L} $ is due to the following  conservation law
\beq\label{conservation}
\gamma_{0}\Bigl( (w,\mathcal{L}  w) + j^2 |\partial_{x}^{-1} w |^2 \Bigr)
= \Re \Bigl( (H_{x}, \mathcal{L}  w) + j^2(H, \partial_{x}^{-1} w )\Bigr)
\eeq
which can be checked by a straightforward computation.
Consequently, we can use  \eqref{dec}, \eqref{decprop}   
and integrate by parts the right-hand side to get
$$
\gamma_{0}\Bigl(\mu\,\alpha(\tau)|\varphi_{-1}|^2+c_{0} |w_{\perp}|^2
+ j^2 |\partial_{x}^{-1} w |^2 \Bigr)
\leq 
C|H|_{2}\, |w|_{1}+ j^2|H| |\partial_{x}^{-1} w|\,.
$$
Therefore, using the inequality
\begin{equation}\label{young}
ab\leq \eps a^2+ \frac{1}{4\eps}\,b^2,\quad \forall\, \eps>0,\quad \forall\, (a,b)\in\mathbb{R}^2,
\end{equation}
with $\eps$ small enough, we can incorporate $|\partial_{x}^{-1} w|$ in the left hand-side and arrive at
\beq\label{wperp}
|w_{\perp}|^2  +  j^2 |\partial_{x}^{-1} w |^2\leq C\Bigl( |\alpha |^2 + 
|H|^2 +  |H|_{2}\, |w|_{1}\Bigr).
\eeq
In what follows $C$ is a large number which may change from lines to lines and depend on $\gamma$  
and $K$ but not on $\tau$. The next step is to estimate $\alpha$ and $\beta$. We use the
decomposition \eqref{dec} and take the scalar product of
\eqref{wL} with $\varphi_{-1}$ and with $\varphi_{0}$ respectively  to get
\begin{eqnarray*}
&& (\gamma_{0}+ i \tau )\alpha   = -(w, \mathcal{L}\partial_{x}(\varphi_{-1})) - j^2(\partial_{x}^{-1}
w , \varphi_{-1}) + (H_{x}, \varphi_{-1}) \\
&&   (\gamma_{0}+ i \tau )\beta  = -(w, \mathcal{L}\partial_{x}(\varphi_{0})) - j^2(\partial_{x}^{-1}
w , \varphi_{0}) + (H_{x}, \varphi_{0})
\end{eqnarray*}
and hence, we can take  the modulus and add the two identities to get
$$ 
(\gamma_{0} + |\tau| ) (  |\alpha | + | \beta | )  \leq C \Bigl( |\alpha | + | \beta| + 
|w_{\perp}| + j^2| \partial_{x}^{-1} w | + | H |_{1} \Bigr).
$$
Next, we multiply  by $|\alpha | + | \beta |$ and use (\ref{young}) to get,
\beq\label{ab}
(\gamma_{0} + |\tau| - C )( |\alpha |^2 + |\beta |^2) \leq C \Bigl(
|w_{\perp}|^2  + j^4| \partial_{x}^{-1} w |^2  + | H |_{1}^2 \Bigr).
\eeq
Note that this last estimate is a good estimate when $\tau$ is large. Next, we can consider  
B\eqref{wperp} +  \eqref{ab} with $B$ a large number to be chosen to get
\begin{multline*}
(B-C)( |w_{\perp}|^2  +  (Bj^2-Cj^4) |\partial_{x}^{-1} w |^2)
+ (\gamma_{0} + |\tau| - C  -BC)( |\alpha |^2 + |\beta |^2)
\\
\leq CB \Bigl(
|H|_{2}\, |w|_{1} + |H|_{1}^2\Bigr).
\end{multline*}
Consequently, we can first choose $B$ sufficiently large (such
 that $B>C, $ and $B>Cj^2$)   and then consider 
$\tau$ sufficiently large (for example $|\tau|\geq 2(C + BC)$) to get the estimate
\beq\label{wL2}
|w|^2 + j^2 |\partial_{x}^{-1} w |^2 \leq C\Bigl(  |H|_{2}\, |w|_{1} + |H|_{1}^2\Bigr), \quad | \tau | \geq M.
\eeq 
To conclude  we just need to  estimate $|\partial_{x} w |$. It suffices to look again at 
\eqref{conservation}. Indeed, we can use that $(w, \mathcal{L}w) = |w_{x}|^2  - \mathcal{O}(1) |w|^2$ 
in \eqref{conservation} to get
\beq\label{wx1}
|w_{x}|^2 + j^2 |\partial_{x}^{-1} w |^2\leq 
C \Bigl( |w|^2 +|H|^2+ |H|_{2} |w|_{1}\Bigr).
\eeq
Consequently, the combination of a sufficiently large constant times \eqref{wL2} and \eqref{wx1} gives
$$ 
|w|_{1}^2 + j^2 |\partial_{x}^{-1} w |^2 
\leq C\Bigl(  |H|_{2}\, |w|_{1} + |H|_{1}^2\Bigr), \quad | \tau | \geq M
$$
and hence by using the inequality (\ref{young}), we  get
\beq\label{H1HF}
|w|_{1}^2 + j^2 |\partial_{x}^{-1} w |^2 \leq C |H|_{2}^2.
\eeq
This proves \eqref{estHF} for $s=1$. Note that moreover \eqref{H1HF} gives
 a control of   $ j^2 |\partial_{x}^{-1}  w |^2$ which is interesting
  when $j \neq 0$.
      
To estimate higher order derivatives, we shall use higher order approximate conservation laws for 
the linearized KdV equation. Namely, we define a self-adjoint  operator
$$ 
\mathcal{L}_{s+1} w = \partial_{x}^{2 s + 2 } w + r_{s+1}(x) \partial_{x}^{2s}w
$$
where $r_{s+1}$  is real valued and will be chosen in order that the following cancellation property occurs : 
\beq\label{Lscom}
\Re  \Bigl((\mathcal{L} w)_{x}, \mathcal{L}_{s+1} w \Bigr) = \mathcal{O}(1) |w|_{s}^2.
\eeq
By making repeated integration by parts, we easily establish that
\begin{eqnarray*}
\Re (\partial_{x}^{2s + 2} w , w_{xxx} ) & = &  (-1)^s\, 
\Re( \partial_{x}^{s+2}w , \partial_{x}\partial_{x}^{s+2}w ) =  0, \\
\Re ( \partial_{x}^{2s + 2} w, Q w_{x})  & =  &  (-1)^{ s + 1} \Re
\Bigl( ( \partial_{x}^{s+1} w , Q \partial_{x}^{s+2} w) + (s+1) ( \partial_{x}^{s + 1} w, 
Q_{x} \,  \partial_{x}^{s+1}w ) + \mathcal{O}(1) | w|_{s}^2\Bigr) \\
&  =  &  (-1)^{ s + 1} \Re\Bigl( (s+ \frac{ 1}{ 2}) ( \partial_{x}^{s+ 1}
w ,  Q_{x }\, \partial_{x}^{s+ 1 } w ) +  \mathcal{O}(1) | w|_{s}^2\Bigr),  \\
\Re ( \partial_{x}^{2s + 2} w, Q_{x} w) & =  &  (-1)^{s + 1} \Re \Bigl(
(\partial_{x}^{s+ 1} w , Q_{x}\, \partial_{x}^{s+1} w \Bigr) + \mathcal{O}(1) |w|_{s}^2 \Bigr),\\
\Re (r_{s+1} \partial_{x}^{2s} w, \partial_{xxx} w )
&  =  & (-1)^{ s - 1 } \Re  \Bigl( \partial_{x}^{ s + 1 } w , \partial_{x}^{s-1}(r_{s+1} w_{xxx}) \Bigr)  \\
& = &    (-1)^{ s - 1 } \Re  \Bigl(  ( \partial_{x}^{s+1}w,  r_{s+1} \partial^{s+2 }_{x} w )  + 
(s-1) (\partial_{x}^{s+1} w , \partial_{x} r_{s+ 1}\partial_{x}^{s+ 1 } w )   \\
& & \mbox{\hspace{2cm}} +  \mathcal{O}(1) |w|_{s}^2 \Bigr)  \\
& =  &    (-1)^{ s - 1 }   \Re  \Bigl(  ( s - \frac{ 3}{ 2})   ( \partial_{x}^{s+1}w, 
\partial_{x} r_{s+1} \partial^{s+1 }_{x} w  ) +     \mathcal{O}(1) |w|_{s}^2 \Bigr) 
\end{eqnarray*}
and that all the other terms which appear in the product $\Re (   ( \mathcal{L} w)_{x}, \mathcal{L}_{s+1}w)$ are 
$ \mathcal{O}(1) |w|_{s}^2$. Consequently, we get
\begin{eqnarray*}
\Re  ( -(\mathcal{L}w)_{x}, \mathcal{L}_{s+1}w) &  = & (-1)^{s+ 1}\Bigl( (s + \frac{ 3}{ 2})  ( \partial_{x}^{s+1}
w ,  Q_{x }\, \partial_{x}^{s+ 1 } w )    \\
& & \mbox{\hspace{1.5cm}} + ( s -  \frac{ 3}{ 2})   ( \partial_{x}^{s+1}w, 
\partial_{x} r_{s+1} \partial_x^{s+1 } w  ) \Bigr) \\
& &     \mbox{\hspace{1.5cm}}    + \mathcal{O}(1) | w |_{s }^2 \\
& = & \mathcal{O}(1) | w |_{s }^2
\end{eqnarray*}
with the choice
$$ 
r_{s+1 } =- \frac{ s + \frac{ 3 }{ 2 }}{ s - \frac{ 3}{ 2} } \, Q .
$$ 
Note that $s$ is an integer so  that $r_{s+1}$ is always well-defined.  
      
Finally, we can  take the scalar product of \eqref{wL} by $(-1)^{s+1}\,\mathcal{L}_{s+1}w$ and  then take
the real part to  get thanks to the above cancellation  property
$$ 
\gamma_{0} |\partial_{x}^{s+ 1 } w |^2 \leq C \Bigl( |w|_{s }^2
 + j^4 |\partial_{x}^{-1} w |^2  + | H |_{s + 2} \, |\partial_{x}^{s+1} w |
+ | H |_{s + 1} \, | w |_{s}
\Bigr) 
$$
since $ \Re  (\partial_{x}^{-1} w,\partial^{2s+ 2}_{x} w) = 0$. We finally obtain
$$
|\partial_{x}^{s+ 1 } w |^2 \leq C \Bigl( |w|_{s }^2 +  j^4 |\partial_{x}^{-1} w |^2
+ | H |_{s + 2}^2 \Bigr) 
$$
thanks to the inequality (\ref{young}) and hence we get \eqref{estHF}  by   induction
 and the control of $ j^2 |\partial_{x}^{-1} w |^2$ given by  \eqref{H1HF}. 
\newline
       
Next, we need to estimate $w$ for $|\tau | \leq M.$ This is the aim of the following lemma.
\begin{lem}\label{BF}
For $|\tau| \leq M$,  we have the estimate
\beq\label{estBF}
|w(\tau)|_{s}^2 \leq C(s,\gamma_{0}, K, M) |H(\tau)|_{s+1}^2.
\eeq
\end{lem}
\subsubsection{Proof of Lemma \ref{BF}}
Note that here we  actually give a proof of the fact that if $\lambda$ is not an eigenvalue 
then $\lambda$ is not in the spectrum.
To prove \eqref{estBF}, we need to treat differently the cases $j=0$ and $j \neq 0$.
  
Let us start with the case $j \neq 0.$ In this case,  we take the derivative of \eqref{wL} to get
\begin{equation}\label{matritza}
(\gamma_{0} + i \tau)w_x - ( \mathcal{L}  w)_{xx} + j^2w = H_{xx} 
\end{equation}
and we introduce $V=(w, w_{x}, w_{xx}, w_{xxx})^t \in \mathbb{C}^4$ and $\mathbb{H}=(0, 0, 0, H_{xx})^t$ 
to rewrite the problem as
\beq\label{abst}
V_{x} = \mathbb{A}(q, x) V + \mathbb{H}
\eeq
where $\mathbb{A}$ is a $4\times 4$ matrix  that one may easily find from the equation (\ref{matritza}) and
the parameter  $q =(\gamma_0+i \tau, j^2)$  is in the compact set $\mathcal{K}$ defined by 
$$
\mathcal{K}= \{(\gamma_0 + i\tau, b), \, |\tau| \leq M, \,  k_{0}^2/ L^2 \leq |b| \leq K^2 k_{0}^2/L^2\}.
$$
Let us denote by $T(q, x,x')$  the fundamental solution of $V_{x} = \mathbb{A}V$ i.e. the solution such that 
$T(q,x',x')=I_{4}$. Next, since $Q(x)$ tends to zero exponentially fast when $x \rightarrow \pm \infty$, 
there exists a matrix $\mathbb{A}_{\infty}(q)$ such that
$$ 
\mathbb{A}(q,x) - \mathbb{A}_\infty(q) = \mathcal{O}(e^{-|x|}), \quad
x\rightarrow \pm \infty.
$$
Moreover   the eigenvalues  of $\mathbb{A}_{\infty}$  are the roots of the  polynomial $P$ defined in 
\eqref{P} below and hence are not purely imaginary.  By classical arguments of ODE (namely the roughness of
exponential dichotomy, see \cite{Coppel} for example), the equation
$V_{x} = \mathbb{A}V$ has an exponential dichotomy on $\rit_{+}$ and $\rit_{-}$, i.e., there exists 
projections $P^+(q, x)$, $P^-(q, x)$ which are smooth in  the parameter with the invariance property 
\begin{equation}\label{invar}
T(q,  x, x') P^\pm(q, x') = P^{\pm}(q , x) T(q, x, x')
\end{equation}
and  such that  there exists $C$ and $\alpha>0$ such that  for every $U \in \mathbb{C}^4$,  
and $q \in \mathcal{K}$, we have 
\begin{eqnarray*}
& &  |T(q, x ,x') P^+(q,  x') U | \leq C e^{-\alpha(x- x')\ } | P^+(q,  x') U|, \, 
x\geq x'\geq 0, \,  \\
& & |T(q, x ,x')( I -  P^+(q, x'))U| \leq C e^{\alpha(x- x')\ } |(I -  P^+(q, x'))U| , \, 
0 \leq x \leq x', \,  \\
& &  |T(q , x ,x') P^-(q,  x') U| \leq C e^{\alpha(x- x')\ }| P^-(q,  x') U| , \, 
x\leq x'\leq 0, \\
& &  |T(q , x ,x') (I - P^-(q,  x')) U| \leq C e^{-\alpha(x- x')\ }| (I - P^-(q,  x') )U| , \, 
0 \geq x\geq x'.
\end{eqnarray*} 
In particular,  note that a solution $T(q, x,0)V^0$ is decaying when $x$ tend to $\pm \infty$ if and only if
$V^0$ belongs to $\mathcal{R}(P^\pm(q, 0))$.
Since by the analysis of \cite{Pego} recalled in section \ref{Appendix} there is no eigenvalue of 
$A_{j}$ (see (\ref{Aj}) for the definition of $A_j$) 
for $q \in  \mathcal{K}$, we have no non trivial solution decaying in both sides and hence we have 
\beq\label{inter}
\mathcal{R}( P^+(q,0) ) \cap \mathcal{R}(P^-(q, 0)  )= \{ 0 \}.
\eeq
Let us choose bases  $(r_{1}^\pm, r_{2}^\pm)$  of $ \mathcal{R}( P^\pm(q,0))$
which depends on the parameters in a smooth way (see \cite{Kato} for example) then we can define
$$ 
M(q ) = (r_{1}^+, r_{2}^+, r_{1}^-, r_{2}^-)
$$ 
and we note that  $M(q)$ is  invertible for $q \in \mathcal{K}$ because of \eqref{inter}.
This allows us to define a new projection $P(q)$  by
$$ 
P(q) = M(q)
\left( 
\begin{array}{cc} 
I_{2} & 0  
\\ 
0 & 0
\end{array} 
\right)
M(q)^{-1}
$$
and next
$$
P(q, x)= T(q,  x, 0) P(q).
$$
The main interest of these definitions is that we have 
 $\mathcal{R}( P(q))= \mathcal{R}( P^+(q,0))$  and $\mathcal{R}( I- P(q))= 
 \mathcal{R}( P^-(q,0))$.  Therefore thanks to (\ref{invar}), we have for every $x$  that  
$\mathcal{R}( P(q,x))= \mathcal{R}( P^+(q,x))$ and  similarly that
$$\mathcal{R}( I - P(q,x)) = \mathcal{R}( P^-(q,x)).$$ 
Consequently, we have the estimates
\begin{eqnarray}
& &    \label{dich+}
|T(q, x ,x') P(q, x')| \leq C e^{-\alpha(x- x')\ } , \, 
x,\, x' \in \mathbb{R}, \, x \geq x',  \,  \, \forall q \in \mathcal{K}, \\
\label{dich-}
& & |T(q, x ,x')(I -  P(q, x'))| \leq C e^{\alpha(x- x')\ } , \, 
x,\, x' \in \mathbb{R}, \, x \leq x',  \,  \, \forall q \in \mathcal{K}.
\end{eqnarray}
By using this property, the unique bounded solution of \eqref{abst} reads by Duhamel formula
$$ 
V(x) = \int_{- \infty}^x T(q ,  x, x')P(q, x') \mathbb{H}(x')\, dx'
- \int_{x}^{+ \infty} T(q, x, x')(I - P(q, x') )\mathbb{H}(x')\, dx'
$$
and hence, we get thanks to \eqref{dich+}, \eqref{dich-} that
$$ 
|V(x)| \leq C \int_{\mathbb{R}} e^{- \alpha| x - x'|} |\mathbb{H}(x') |\, dx'
$$
which yields by standard convolution estimates 
$$ 
|V| \leq C |H|.
$$
The estimates of high order derivatives is very easy, it suffices to write
$$ 
\partial_{x}^{s+1} V = \mathbb{A} \partial_{x}^sV + [ \partial_{x}^s, \mathbb{A}]V +\partial_{x}^s \mathbb{H},
$$
and to write Duhamel formula considering $[ \partial_{x}^s, \mathbb{A}]V$ as part of the source term.
\newline
        
It remains the case $j=0$. In this case, we do not take the derivative of \eqref{wL}, we directly 
define $W=(w, w_{x}, w_{xx})$ and we rewrite \eqref{wL} under the form
$$ 
W_{x} = \mathbb{B}(\lambda, x) W+ \tilde{\mathbb{H}}.
$$
Then the proof of the estimate follows the same line, we find that $\mathbb{B}_{\infty}$ has no 
eigenvalue on the imaginary axis. This yields that there is an exponential dichotomy on $\mathbb{R}_{+}$
and $\mathbb{R}_{-}$ for this system. Next since, the spectrum of the linearized KdV equation about the 
soliton is on the imaginary axis, we get that the system has an exponential
dichotomy on the real line. We do not detail more since the proof is similar to the previous case.
\subsubsection{End of the proof of Theorem \ref{resolvant} }
To get \eqref{estresolvant}, it suffices to combine Lemma  \ref{HF} and Lemma \ref{BF}.
\subsubsection{End of the proof of Theorem \ref{theolinsource}}
By using Theorem \ref{resolvant} and Bessel-Parseval identity, we get that for every $T>0$, 
\begin{eqnarray*}
& & \int_{0}^T e^{-2 \gamma_{0} t } |v(t)|_{s}^2\, dt
\leq \int_{0}^{+ \infty} e^{-2 \gamma_{0} t } |\tilde{v}(t)|_{s}^2\, dt
= \int_{\mathbb{R}} |w(\tau)|_{s}^2 \, d\tau \\
& & \leq C \int_{\mathbb{R}} |H(\tau)|_{s+1}^2 \, d\tau
= \int_{0}^{T} e^{-2 \gamma_{0} t } |F_j(t)|_{s+1}^2\, dt
\end{eqnarray*}
and finally thanks to \eqref{F}, we get
\beq\label{estsource}
\int_{0}^T e^{-2 \gamma_{0} t } |v(t)|_{s}^2\, dt\leq C \int_{0}^T e^{2( \gamma - \gamma_{0})t}\, dt
\leq C e^{ 2( \gamma - \gamma_{0})T}
\eeq
since $\gamma_{0}$ was fixed  such that $\gamma> \gamma_{0}$.
To finish the proof, we notice that the energy estimate for the equation \eqref{ulink} gives
$$ 
\frac{ d}{ dt} |v(t)|_{s}^2 \leq C \Bigl( |v(t)|_{s}^2 + |F_j(t)|_{s+1}^2\Bigr).
$$
Consequently, we can multiply the last estimate by $e^{-2\gamma_{0}t}$ and use \eqref{Fbis} to get
$$
\frac{ d}{  dt} \Bigl( e^{-2 \gamma_{0} t } |v(t)|_{s}^2\Bigr)
\leq C \Bigl( e^{- 2 \gamma_{0}t } |v(t)|_{s}^2  + e^{2 (\gamma - \gamma_{0})t}\Bigr).
$$
Next, we integrate in time and use \eqref{estsource} and again  the fact that $\gamma>\gamma_{0}$, this yields
$$ e^{- 2 \gamma_{0} t }  |v(t)|_{s}^2 \leq C e^{2 (\gamma - \gamma_{0}) t }. $$
This ends the proof of Theorem~\ref{theolinsource} .
\subsubsection{Proof of Proposition~\ref{propuk}}
By induction, it suffices to use Theorem~\ref{theolinsource} and 
the fact that $H^s(\mathbb{R})$ is an algebra for $s\geq 1$.

\subsection{Nonlinear instability: end of the proof of Theorem~\ref{maintheo}}     
Of course, we only need to prove the statement for $\delta$ small enough.
Let us define $w$ by setting $v=u^{ap}+w$, where $u^{ap}$ is defined by (\ref{u_app}). 
Therefore we have that the solution $u^\delta$ may be decomposed as follows
$$
u^\delta = Q+u^{ap}+w\,.
$$
If we set
$$
F\equiv (\partial_t+A)u^{ap} + u^{ap} u^{ap}_{x}\,, 
$$
where $A$ is defined in (\ref{eqlin}), then thanks to Proposition~\ref{propuk},
$$
\|F(t,\cdot)\|_{L^2(\rit\times \tit_{L})}\leq C_{M,s}\delta^{M+2}e^{(M+2)\sigma_0 t}.
$$
We have that $w$ solves the problem
\begin{equation}\label{ww}
(\partial_t+A)w+\partial_{x}(u^{ap}w)+ww_x+F=0\,, \quad w_{/t=0}= 0 \,.
\end{equation}
We now estimate the solution of (\ref{ww}). Using that
$$
|\int_{ \rit\times \tit_{L}}F w|\leq \|F(t,\cdot)\|_{L^2(\rit\times \tit_{L})}^{2}+
\|w(t,\cdot)\|_{L^2(\rit\times \tit_{L})}^2
$$
multiplying (\ref{ww}) by $w$ and integrating $\rit\times \tit_{L}$, 
we get after several integrations by parts
\begin{equation}\label{energia}
\frac{d}{dt}\|w(t,\cdot)\|_{L^2}^{2}\leq
\Big(\|Q'\|_{L^{\infty}}+\|\partial_{x}u_{ap}(t,\cdot)\|_{L^\infty}+1\Big)\|w(t,\cdot)\|_{L^2}^{2}
+\|F(t,\cdot)\|_{L^2}^{2}\,.
\end{equation}
Observe that
$$
\|\partial_{x}u_{ap}(t,\cdot)\|_{L^\infty(\rit\times \tit_{L})}\leq
\sum_{k=0}^{M}C_{k,s}\delta^{k+1}e^{(k+1)\sigma_0 t}\,.
$$
Next, we set
$$
T^\delta\equiv \frac{\log(\kappa/\delta)}{\sigma_0}\,,
$$
where $\kappa\in]0,1]$ is small enough to be chosen after the several restrictions we will impose in 
the next lines. The number $T^\delta$ represents the time when the instability occurs.
Coming back to (\ref{energia}),
we observe that there exists a constant $\Lambda_{M,s}$ 
depending on $s$ and $M$ but independent of $\kappa$ and $t$ 
and an absolute constant $C$ ($C$ is essentially $\|Q'\|_{L^{\infty}}$)  such that for 
$0\leq t\leq T^\delta$,
$$
\frac{d}{dt}\|w(t,\cdot)\|_{L^2}^{2}\leq
(C+\kappa\Lambda_{M,s})\|w(t,\cdot)\|_{L^2}^{2}
+ C_{M,s}\delta^{2(M+2)}e^{2(M+2)\sigma_0 t}\,.
$$
Therefore
\begin{equation}\label{veneta}
\frac{d}{dt}\Big(e^{-(\kappa\Lambda_{M,s}+C) t}\|w(t,\cdot)\|_{L^2}^{2}\Big)
\leq 
C_{M,s}\delta^{2(M+2)}e^{2(M+2)\sigma_0 t-\kappa\Lambda_{M,s} t-Ct}\,,\quad t\in [0,T^\delta]\,.
\end{equation}
Now we choose $M$ large enough and $\kappa$ small enough so that 
$$
2(M+2)\sigma_0-\kappa\Lambda_{M,s}-C>0.
$$
At this place we fix the value of $M$ (and of $s$, for exemple $s=M+1$) 
while we will make two more restrictions on $\kappa$.
Since $w$ vanishes for $t=0$ an integration of (\ref{veneta}) yields
$$
\|w(t,\cdot)\|_{L^2(\rit\times \tit_{L})}\leq 
C_{M,s}\delta^{M+2}e^{(M+2)\sigma_0 t}\,,\quad t\in [0,T^\delta]\,.
$$
Therefore
\begin{equation}\label{reste}
\|w(T^\delta,\cdot)\|_{L^2(\rit\times \tit_{L})}\leq C_{M,s}\kappa^{M+2}\,.
\end{equation}
Let us denote by $\Pi$ the projection on the nonzero modes in $y$ i.e.
$$(\Pi v)(x,y) \equiv  v(x,y) - { 1 \over 2 \pi L } \int_{0}^{2\pi L} v(x,y)\, dy.$$
 Then for every $a\in \mathbb{R}$ one has
$\Pi(Q(x-a))=0$. On the other hand the first term of $u^{ap}$ satisfies $\Pi(u^{0})=u^{0}$ and 
therefore
$$
\|\Pi(u^{ap}(t,\cdot))\|_{L^2}
\geq
c_s\delta e^{\sigma_0 t}-\sum_{k=1}^{M}\delta^{k+1}\|\Pi(u^{k})\|_{L^2}
\geq
c_s\delta e^{\sigma_0 t}-\sum_{k=1}^{M}C_{k,s}\delta^{k+1}e^{(k+1)\sigma_0 t}\,,
$$
where $c_s$ is the $H^s(\mathbb{R}\times \mathbb{T}_{L})$ norm of $u^0$.
Therefore for $\kappa$ small enough one has 
\begin{equation}\label{parvo}
\|\Pi(u^{ap}(T^\delta,\cdot))\|_{L^2(\rit\times \tit_{L})}
\geq
\frac{c_s \kappa}{2}\, .
\end{equation}
Using (\ref{reste}) and (\ref{parvo}), we may write that for every $a\in\mathbb{R}$,
\begin{multline*}
\|u^{\delta}(T^\delta,\cdot)-Q(\cdot-a)\|_{L^2}
\geq 
\|\Pi(u^{\delta}(T^\delta,\cdot)-Q(\cdot-a))\|_{L^2}
\\
=\|\Pi(u^{\delta}(T^\delta,\cdot)-Q(\cdot))\|_{L^2}
=\|\Pi(u^{ap}(T^\delta,\cdot)+w(T^\delta,\cdot))\|_{L^2}
\\
\geq \frac{c_s \kappa}{2}-\|\Pi(w(T^\delta,\cdot))\|_{L^2}
\geq \frac{c_s \kappa}{2}-\|w(T^\delta,\cdot)\|_{L^2}
\geq
\frac{c_s \kappa}{2}-C_{M,s}\kappa^{M+2}\,.
\end{multline*}
A final restriction on $\kappa$ may insure that the right hand-side of the last inequality 
is bounded from below by a fixed positive constant $\eta$ depending only on $s$
(in particular $\eta$ is independent of $\delta$).
This completes the proof of Theorem~\ref{maintheo}.
\begin{rem}
Let us observe that the analysis in the proof of Theorem~\ref{maintheo} is quite different from the 
high frequency instabilities studied in \cite{KT}.
In \cite{KT},  the approximated solution is a high frequency linear wave with modified speed, perturbed
by a low frequency wave. 
In Theorem~\ref{maintheo},  the approximated solution is a low frequency object modelled on the profile
$u^0$. 
\end{rem}
\section{Proof of Theorem~\ref{theoschro}}
The proof follows exactly the same lines as the proof on Theorem~\ref{maintheo}
and thus we shall only sketch it. We again  look for $u^\delta$ under the form
$u^\delta = Q + u^{ap} + w$.  At first, we  need to find  a most unstable
eigenmode for the linearized equation to begin the construction
 of $ u^{ap}$. The linearized equation about $Q$
 reads
 $$i u_{t} + A u = 0, \quad A u = \Delta u - u + 2 u Q^2 + \overline{u} Q^2.$$
 It is more convenient to introduce $ U= (\Re u, \Im u)^t$ and to rewrite
  the equation as the system :
  \beq
  U_{t} +  \left( \begin{array}{cc} 0  & - \mathcal{ L}^{-} \\   \mathcal{L}^{+} 
  & 0 \end{array}\right) U = 0, 
  \eeq 
$$ 
\mathcal{L}^{-} u =  - \Delta u  + u - Q^2 u , \quad \mathcal{L}^{+} u = -\Delta u  +u - 3 Q^2 u.
$$ 
We seek unstable  eigenmodes under the form 
\begin{equation}\label{form}
\Phi_{\sigma, k }(t,x,y ) = 
e^{\sigma t}
e^{ i  k y \over L} V(x)+
e^{\overline{\sigma} t}
e^{{- i  k y \over L}}\, \overline{V}(x), \quad  \Re \sigma >0,
\end{equation}
where $V(x)\in {\mathbb C}^2$ so that we have to solve
\beq\label{eigenmodes}
\sigma V +\left( \begin{array}{cc} 0  & -  L^{-} - {k^2
 \over L^2}  \\   L^{+} + {k^2 \over L^2} & 0 \end{array}\right) V = 0
\eeq
where
$$ 
L^{-} u =  - u_{xx}  + u   - Q^2 u , \quad L^{+} u = - u_{xx}  + u   - 3 Q^2 u.
$$ 
We set $\eps = {k  \over L}$ and we look for nontrivial solutions of
 \eqref{eigenmodes} with $\Re \sigma >0$ for $\eps>0.$
The first result we shall use is that
\begin{lem}
\label{specschro}
For $\eps>0$, there is at  most one unstable eigenmode and
 there exists $\eps_{0}$  such that for $0<\eps \leq \eps_{0}$, there 
is exactly one unstable eigenmode. 
\end{lem}
In the reference \cite{JR}, it is claimed that the result of this lemma is due to
Zakharov-Rubenchik. Unfortunately, we were not able to find a copy of the paper by Zakharov-Rubenchik
as this paper is quoted in \cite{JR}. We give a proof of this Lemma in the appendix.
\\

Now,   thanks to  Lemma \ref{specschro}, for $ k=1$ and $L$ sufficiently large
there exists an unstable eigenmode. We  now consider $L$ as fixed.  For every
$k$, we have by Lemma~\ref{specschro} that  there exists at most one
$\sigma(k)$ such that $ \Re \sigma (k)>0$ and (\ref{eigenmodes}) has a
solution in $L^2({\mathbb R};{\mathbb C}^2)$ with $\sigma=\sigma(k)$.
Moreover we can easily get that the solutions of (\ref{eigenmodes}) satisfy
the conservation law 
$$  
\Re \sigma \Bigl( (L^{+ }V_{1}, V_{1}) + (L^-V_{2}, V_{2} ) + {k^2 \over L^2} |V|^2
 \Bigr)  =0\,.
$$
Therefore for large $k$ (depending only on $Q$) there is no nontrivial solution of (\ref{eigenmodes}) with $\Re \sigma >0$.
Consequently, we can choose an eigenmode $\Phi_{\sigma, k }$ under the form \eqref{form} such that
  $$  \Re \sigma = \sup \Bigl\{ \Re  \sigma(k)  \Bigr\} := \sigma_{0}$$
and we set $u^0= (\Phi_{\sigma,k})_{1} + i (\Phi_{\sigma,k})_{2}$.
Observe that thanks to (\ref{eigenmodes}) we have $(i\partial_t+A)u^0=0$.
The next step towards the  proof of  Theorem~\ref{theoschro} is the construction on an high order unstable solution.
We use the same method as previously, we use the same spaces $V_{K}^s$ and we build 
an approximate solution under the form \eqref{u_app}. For $ 1 \leq k \leq M+1$, we need to solve
\beq\label{ukschro}
i \partial_{t }u^k  + A u^k = -   \sum_{j+l=k-1}( 2Q u^j \overline{u}^l + Q
u^j u^l) - \sum_{j+l+m=k-2} u^j \overline{u}^l u^m, \quad (u^k)_{/t=0}=0 
\eeq
where the last sum is zero for $k=1$.  We have the estimates :
\begin{prop}\label{propukschro}
Let $u^k$ the solution of \eqref{ukschro}, we have the estimate
$$ |u^k(t) |_{V^s_{k+1}} \leq C_{k,s} e^{ (k+1) \sigma_{0} t }, \, \forall t \geq 0.$$
\end{prop}
Note that here we do not loose regularity at each step because the nonlinear term does not involve derivatives.
To prove Proposition~\ref{propukschro},  
we need to prove the  equivalent of Theorem~\ref{theolinsource}. By using Laplace transform,
we can still reduce the problem to the proof of  a resolvent estimate as in Theorem~\ref{resolvant}.
The proof of the low frequencies estimates rely on 
the same ODE argument and we shall not detail it. We shall just explain how to get
 the high  frequencies estimates.  As in Lemma~\ref{specschro}, it is more convenient to work on 
the  system form of the problem, and thus we consider the equation
\beq\label{resolschro}
(\gamma_{0} + i \tau) W + \left( \begin{array}{cc} 0  & -  L^{-}-\frac{k^2}{L^2} \\   
L^{+}+\frac{k^2}{L^2} & 0 \end{array}\right) W = H
\eeq
and we want to prove that $W(\tau)$  satisfies the estimate
\beq\label{HFschro}
|W(\tau)|_{s}^2 \leq C(s,\gamma_{0}, K) |H(\tau)|_{s}^2
\eeq
for $\gamma_0>\sigma_0$, $ |\tau |\geq M\gg 1$ and $ s\geq 1$.  We first give the proof for $s=1$. 
The conservation law reads  for $W=(w_{1}, w_{2})$ 
\beq\label{conschro} 
\gamma_{0} \Bigl( (L^{+ }w_{1}, w_{1}) + (L^-w_{2}, w_{2} ) + {k^2 \over L^2} |W|^2
 \Bigr) = \Re\Big((H_{1}, L^+ w_{1}) + (H_{2}, L^- w_{2})\Big). 
\eeq
At this stage, we shall use the description of the spectrum of $L^{\pm}$
recalled in the appendix of this paper. We can write
$$
w_{2} = \alpha Q + w_{2}^\perp, \quad  (L^-w_{2}^\perp, w_{2}^\perp) 
\geq c_{0} | w_{2}^\perp |^2.
$$
Similarly, we can write
$$ 
w_{1} = \beta \varphi_{-1} + \gamma Q_{x} + w_{1}^{\perp}, \quad (L^+w_{1}^\perp, w_{1}^\perp) 
\geq c_{0} | w_{1}^\perp |^2 
$$
($\varphi_{-1}\equiv Q^2$).
Setting $W^\perp = (w_{1}^\perp, w_{2}^\perp)^t$ and $W_{F}=(\alpha, \beta, \gamma )^t\in \mathbb{C}^3$, 
we get from \eqref{conschro}
\beq\label{schro1}
\gamma_{0} | W^{\perp}(\tau) |^2 \leq C \Bigl( |H|_{1} \, |W|_{1} + |W_F|^2\Bigr).
\eeq
Next, we can take the projection of the equation on the finite dimensional subspace generated by  
$(0, Q)$, $(Q_{x}, 0)$, $(\varphi_{-1}, 0)$ to get
\beq\label{schro2}
(\gamma_{0} + |\tau| - C ) |W_{F}|^2  \leq  C(K) \Bigl( |W^\perp |^2 +  |H|^2\Bigr).
\eeq
As for the KP-I equation, a suitable combination of \eqref{conschro}, 
\eqref{schro1}, \eqref{schro2}  with the use of \eqref{young} gives \eqref{HFschro} for $s=1$ for $|\tau|$ large
enough. To get higher order derivatives, we use approximate higher order conservation laws. 
Namely, we choose $L^+_{s+1}$ and $L^-_{s-1}$ under the form
$$  
L^-_{s+1} w= \partial_{x}^{2(s+1)}w + r_{s+1}^-\partial_{x}^{2s}w \, \quad L^+_{s+1} w =  
\partial_{x}^{2(s+1)}w + r_{s+1}^+\partial_{x}^{2s}w 
$$ 
such that the cancellation 
$$  
- \Re\Big(L^-w_{2}  + {k^2 \over L^2} w_{2} ,  L^+_{s+1} w_{1} \Big) +  \Re\Big(L^+w_{1}  + 
{k^2 \over L^2} w_{1} ,  L^-_{s+1} w_{2} \Big) = \mathcal{O}(1)( |W|_{s}^2 + |W|_{s} \, 
|\partial_{x}^{s+1}W|)
$$
occurs. To perform this cancellation, it suffices to choose $r_{s+1}^{\pm}$
real valued and such that
$$ 
r_{s+1}^+ - r_{s+1}^- = 2Q^2.
$$
Using this approximate conservation law, we get
$$
\gamma_{0} |\partial_{x}^{s+1}W|^2 \leq C \Bigl( |W|_{s}^2 + 
  |W|_{s} \, 
|W|_{s+1} + |H|_{s+1} \, |\partial_{x}^{s+1}W|
+
|H|_{s} \, |W|_s
\Bigr)$$
and we  conclude thanks to \eqref{young} via an induction argument.
\\

To end the proof of Theorem \ref{theoschro}, we seek for a solution of \eqref{NLS} 
under the form $ u^\delta = Q + u^{ap}+w$, with $w_{/t=0} = 0$ so that $w$ solves the equation 
\beq
\label{wschro}
i w_{t} + A w + 2 |u^{ap}|^2 w +  (u^{ap})^2 \overline{w} + \mathcal{N}(u^{ap}, w) + |w|^2 w = F
\eeq
with
$$ 
\|F\|_{H^s(\rit \times \tit_{L})} \leq C_{M,s} \delta^{M+2} e^{(M+2) \sigma_{0} t},
$$
and  the bilinear  term satisfies
$$ 
\|\mathcal{N}(u^{ap}, w)\|_{H^s(\rit \times \tit_{L})} \leq  C |u^{ap}|_{W^{s,\infty}} \| w\|_{s}^2.
$$
Since here we do not have a global existence result available, we shall first prove that this 
last equation has a smooth solution $w$ which remains defined on a time scale sufficiently 
long to see the instability.
     
A classical existence result   for this equation based on Duhamel formula and Sobolev embedding 
gives that there exists a local  solution $w \in \mathcal{C}([0,T], H^s)$  for $s>1$. Moreover,
we can define a maximum time $T^*$ such that
$$ 
T^*= \sup \{ T, \, \forall t \in [0,T], \|w(t)\|_{H^s} \leq 1 \}.
$$
The $H^s$ energy estimate for \eqref{wschro}  gives for $t \in [0, T^*)$  that
$$ 
{d \over dt} \| w(t) \|_{H^s}^2 \leq C ( 1 + |u^{ap}|_{W^{s,\infty} } ) \|w\|_{H^s}
+ C_{M,s}\delta^{2(M+2)}e^{2(M+2)\sigma_0 t}
$$
where $C$ is an absolute constant (which depends on $Q$). 
Consequently for 
$$
t \leq \mbox{ Min }( T^\delta : = {\log(\kappa /\delta) \over \sigma_{0} }, T^*)\, ,
$$
we get
$$
{d \over dt} \| w(t) \|_{H^s}^2 \leq  ( C + \kappa \Lambda_{M,s} ) \|w\|_{H^s}
+ C_{M,s}\delta^{2(M+2)}e^{2(M+2)\sigma_0 t}
$$
and hence by the choice
$$ 
2 (M+ 2 ) \sigma_{0} - \kappa \Lambda_{M,s} - C >0,
$$
we get that
\beq\label{wschrofin}
\|w(t)\|_{H^s(\rit\times \tit_{L})}\leq C_{M,s}\kappa^{M+2}, \, \quad t \leq   
\mbox{ Min }( T^\delta , T^*).\eeq
In particular for  $\kappa$ sufficiently small, we get that
$$ \|w(t)\|_{H^s(\rit\times \tit_{L})}\leq {1 \over 2}, \quad t \leq \mbox{ Min }( T^\delta, T^*).$$
By definition of $T^*$, this proves that $T^* \geq T^\delta$ so that the time of existence 
of a smooth solution is in any case large enough to see an instability. 
The end of the proof follows the same lines as previously, using again the projection $\Pi$ on nonzero modes
in $y$,  we write for every $a\in\mathbb{R}, \, \gamma \in \rit$,
\begin{multline*}
\|u^{\delta}(T^\delta,\cdot)-e^{i \gamma }Q(\cdot-a)\|_{L^2}
\geq 
\|\Pi(u^{\delta}(T^\delta,\cdot)-e^{i \gamma }Q(\cdot-a))\|_{L^2}
=|| \Pi ( u^{ap}(T^\delta, \cdot) + w(T^{\delta}, \cdot)) ||_{L^2}
\\
\geq \frac{c_s \kappa}{2}-\|\Pi(w(T^\delta,\cdot))\|_{H^s}
\geq \frac{c_s \kappa}{2}-\|w(T^\delta,\cdot)\|_{H^s}
\geq
\frac{c_s \kappa}{2}-C_{M,s}\kappa^{M+2}\,.
\end{multline*}
where we have used \eqref{wschrofin} in the last inequality. A final restriction 
of $\kappa$ gives the instability result.

\section{Appendix} 
\subsection{ Proof of Theorem \ref{Pego}}
\label{Appendix}      
In order to have the same  equations as in \cite{Pego}, we look for solutions of \eqref{eqlin} under the form
$$ 
u(t,x,y)= e^{\frac{\lambda t}{ 2}} e^{\frac{iky}{ L}} U\Bigl(\frac{x}{2}\Bigr)
$$
with $U \in L^2$, $\Re \lambda >0$ and $k\neq0.$ Note that this last condition is natural since for $k=0$,  
we cannot find instability since the KdV soliton is stable in the KdV equation. We get for $U$ the equation
\beq\label{U}
4 \lambda U_{z} + 4 ( \Phi\, U)_{zz} + U_{zzzz}-4U_{zz} +  3 \eta^2 U = 0
\eeq
where we have set
\beq\label{eta}
3 \eta^2 = \frac{ 16 k^2}{ L^2}
\eeq 
and $\Phi = 3 \mbox{ sech}^2 z$. Since $\Phi$ and its derivatives tend to zero exponentially fast when 
$z\rightarrow \pm \infty$,
the solutions of \eqref{U}  have the same behaviour as the solutions of
$$ 
4 \lambda U_{z}  + U_{zzzz}-4U_{zz} +  3 \eta^2 U = 0
$$
when $z\rightarrow \pm \infty$. The characteristic values $\mu$ of this linear equation are the  
roots of the polynomial $P$ defined by 
\beq\label{P}
P(\mu ) = \mu^4 - 4 \mu^2 + 4 \lambda \mu + 3 \eta^2.
\eeq
Consequently for $\eta \neq 0$ and $\gamma = \Re \lambda >0$, $\mu \notin i \mathbb{R}$. Indeed, 
if $\mu = i\xi \in i\mathbb{R}$, then $\xi$ should solve 
$$
\xi^4+4\xi^2+4\lambda\xi i+3\eta^2=0
$$
which cannot have a real root $\xi$ for $\eta\neq 0$ and $\Re\lambda \neq 0$.
A consequence of this is that  the number of roots $\mu$ of positive real part  of 
$P$ is independent of the parameters. Since the limit $\eta\rightarrow + \infty$ gives
$$ 
\mu = 3^{ \frac{1}{ 4}}\, \omega \sqrt{\eta} + \mathcal{O}(1), \quad\omega^4 = -1
$$
we finally get that $P$ has two roots of positive real parts and two roots of negative real parts. 
This proves that   the solutions of \eqref{U} either tends to  zero or blows-up exponentially fast
when $z \rightarrow \pm \infty$. Moreover, the stable manifold and the  unstable manifold have the same  
dimension $2$.  Finally,  there will be a nontrivial bounded solution of \eqref{U} if and only if
$U$ belongs simultaneously to the stable and the unstable manifold.
              
In our case, this condition  can be computed explicitly. Indeed, we notice that for $\gamma>0$, 
$\eta \neq 0$ there is a bounded solution of \eqref{U} if and only if $U= g_{zz}$ with $g$ bounded 
which solves
\beq\label{g}
g_{zzzz} + 4\,\Phi\,g_{zz}  + 4 \,\lambda g_{z}-4g_{zz} + 3 \,\eta^2 g = 0
\eeq
Note that the asymptotic behaviour of the solutions of this equation is also determined by the 
characteristic values given by the roots of $P$ so that this equation also has stable and unstable manifolds
of dimension 2. Moreover, if $\mu$ is a root of $P$, then
\beq
g_{\mu}(z) = e^{\mu z}\Bigl( \mu^3 + 2 \mu + \lambda - 3\mu^2 \,\mbox{tanh } z\Bigr)
\eeq
is a solution of \eqref{g}. In particular, if $ \Re \mu>0$, then $g_{\mu}$ is in the unstable manifold.  
Moreover, when  $P$ has two simple roots $\mu_{1}, $ $\mu_{2}$ of positive real parts, then one can 
prove (see \cite{Pego} for details) that $g_{\mu_{1}}, $ $  g_{\mu_{2}}$ are linearly independent
so that they constitute a basis of the unstable manifold. Consequently, any bounded solution of 
\eqref{g} must be a linear combination of  $g_{\mu_{1}}, $ and  $  g_{\mu_{2}}$.
   
Now, let us define
$$ 
C_{+}(\mu) = \lim_{z \rightarrow + \infty} e^{-\mu z} g_\mu =  \mu^3 + 2 \mu + \lambda - 3\mu^2.
$$
Then,  if $C_{+}(\mu_{i}) \neq 0$, $i=1, \, 2,$  we cannot have nontrivial
  solutions which tend to zero when 
$z \rightarrow + \infty$. Consequently, this proves that when
the positive real part roots of $P$ are simple, then  a necessary condition to have
bounded solutions of \eqref{g} is that $C_{+}(\mu)$ = 0 for some root $\mu$ of $P$ of positive real part. 
In the case where $\mu$ is a double root, then  one can check
that the same condition holds. Indeed it suffices to take $g_{\mu}$ and $\partial_{\mu}g$ as 
a basis of the unstable manifold (again, we refer to \cite{Pego} for details).
                              
It remains to study the equation $C_{+}(\mu)$ = 0 with $\mu$ a root of $P$ of positive real part. This yields
the system of algebraic equation
\beq\label{alg}
P(\mu)= 0, \quad \mu^3 + 2 \mu + \lambda - 3 \mu^2 = 0,
\eeq
with the constraint $\Re \mu >0$. The elimination of $\lambda$ between the two algebraic equations gives
\beq\label{alg2}
\lambda = - \mu (\mu -1) (\mu - 2), \quad \eta^2 = \mu^2(\mu - 2 ) ^2.
\eeq
The analysis of this system gives that  there is a solution with $\Re \lambda >0$, $\Re\mu >0$, if and only if
given $\mu \in (0, 2)$,  $\eta$ and $\lambda$ are given by 
$$ 
\eta = \mu(2 - \mu), \quad \lambda = -\mu(\mu-1)(\mu -2).
$$
Finally, we notice   that when $C^+(\mu)=0$, we have
$$ 
g_{\mu}(z) = 3 \mu^2 e^{\mu z}( 1- \mbox{tanh }z ) = \mathcal{O}( e^{- (2 - \mu) z})
$$
and hence $\lim_{z\rightarrow + \infty} g_{\mu} = 0$ since $2 - \mu >0$.
This proves that $C_{+}(\mu)=0$ with $\mu$ a root of $P$ of positive real part is also a sufficient 
condition to have a bounded solution on $\mathbb{R}$. This ends the proof.
\subsection{Proof of Lemma \ref{specschro}}
Set   $ V(x)= (u(x), v(x))^t$ with $u,v$ real valued functions.
Then \eqref{eigenmodes} implies that
\begin{equation}\label{eig-bis}
 L^{+}u+\eps^2 u=-\sigma v,\quad L^{-}v+\eps^2 v=\sigma u\,.
\end{equation}
Observe that if $(u,v)$ is a solution of (\ref{eig-bis}) corresponding to a complex number $\sigma$
then $(u,-v)$ is a solution of (\ref{eig-bis}) corresponding to $-\sigma$.
The operators $L^+$ and $L^-$ have classical self adjoint realizations on $L^2(\mathbb{R})$ and
their spectrum are well-known (see e.g. \cite{T,W}). The operator $L^+$ has exactly two simple eigenvalues
$-3$ and $0$ with corresponding eigenfunctions $Q^2$ and $Q'$.  The continuous spectrum of   
$L^+$ is $[1,\infty[$. The operator $L^-$ has only the simple eigenvalue $0$ with corresponding
eigenfunction $Q$ and the continuous spectrum of $L^-$ is $[1,\infty[$.
Observe that (\ref{eig-bis}) may be written as
\begin{equation}\label{op}
 \mathcal{L} \left( \begin{array}{ll} u \\ v \end{array} \right) : =  \left( \begin{array}{cc} 0  & - 1 \\   1 
  & 0 \end{array}\right)
\left( \begin{array}{cc} L^+ +\eps^2  & 0 \\   0 
  & L^- +\eps^2 \end{array}\right)
\left( \begin{array}{c} u   \\   v 
 \end{array}\right)
=
-\sigma\left( \begin{array}{c} u   \\   v 
 \end{array}\right).
\end{equation}
Thanks to the above discussion on the spectrum of $L^+$ and $L^-$, we obtain that
$$
\left( \begin{array}{cc} L^+ +\eps^2  & 0 \\   0 
  & L^- +\eps^2 \end{array}\right)
$$
has at most one negative eigenvalue which should be simple. Therefore, thanks to
\cite[Theorem~3.1]{PW}, there cannot be more than one unstable mode. 

For $\eps\ll 1$, the bifurcation of the eigenvalue zero in the case $\eps=0$
can be explicitly computed. Note that zero is an isolated eigenvalue so that
 we can use perturbation methods as in finite dimension (see
  \cite{Kato} Theorem 1.8, Chapter 7). 
In the case $\eps=0$, we have that zero is an eigenvalue of multiplicity $4$ for the linear map 
introduced in the left hand-side of (\ref{op}) (see \cite{W}).
The generalized eigenspace splits into two two dimensional invariant sub-spaces corresponding to the 
eigenvectors $(u,v)=(Q',0)$ and $(u,v)=(0,Q)$ respectively.
 As generalized eigenvectors, we can take
 $ { 1 \over 2} (Q+ xQ_{x}, 0)$ and $ (0, { 1 \over 2} xQ)$
  which verify 
 $$ \mathcal{L} \left( \begin{array}{cc} {1 \over 2} (Q + x Q_{x})
   \\ 0 \end{array} \right) = - \left( \begin{array}{cc} 0 \\ Q
   \end{array} \right), \quad
   \mathcal{L} \left( \begin{array}{cc} 0
   \\ { 1 \over 2 } x Q \end{array} \right) =  \left( \begin{array}{cc} Q_{x} \\ 0 
   \end{array} \right).$$
 
Thanks to the analytic dependence in $\eps$ (see \cite{Kato}), we look for a 
$\sigma$ in (\ref{eig-bis}) of the form
$
\sigma=\omega_1\eps+\omega_2\eps^2+\cdots 
$
with $\Re(\omega_1)>0$ which corresponds to an unstable mode.
We will see below that the invariant subspace corresponding to 
$(u,v)=(Q',0)$ splits to two one-dimensional invariant spaces corresponding to eigenvalues
with $\omega_1$ purely imaginary and, what is of importance for our purposes, 
the invariant subspace corresponding to 
$(u,v)=(0,Q)$ splits to two one-dimensional invariant spaces corresponding to eigenvalues
with positive and negative $\omega_1$. The eigenvector corresponding to a positive $\omega_1$
provides the unstable eigenmode.
Assume that $u$ and $v$ are expanded as
$$
u=u_0+u_1\eps+u_2\eps^2+\cdots,\quad
v=v_0+v_1\eps+v_2\eps^2+\cdots\,.
$$
Then $(u_0,v_0)$ satisfy $L^+(u_0)=L^- (v_0)=0$.
Thus there exist two numbers $\alpha_0$ and $\beta_0$ such that
$u_0=\alpha_0 Q'$ and $v_0=\beta_0 Q$.
Then $(u_1,v_1)$ are solutions of $L^+(u_1)=-\omega_1 \beta_0 Q$, 
$L^-(v_1)=\omega_1 \alpha_0 Q'$. 
Therefore there exists two numbers $\alpha_1$ and $\beta_1$ such that
$$
u_1(x)=\frac{\omega_1\beta_0}{2}(xQ'(x)+Q(x))+\alpha_1 Q'(x),\quad 
v_1(x)=
-\frac{\omega_1\alpha_0}{2}(xQ(x))
+\beta_1 Q(x)\,.
$$
Next, $(u_2,v_2)$ are solutions of 
\begin{multline}\label{v2}
L^+(u_2)=
-\alpha_0 Q'
-\omega_1\Big(-\frac{\omega_1\alpha_0}{2}(xQ)
+\beta_1 Q\Big)
-\omega_2 \beta_0 Q,
\\
L^-(v_2)=-\beta_0 Q
+\omega_1\big(\frac{\omega_1\beta_0}{2}(xQ'+Q)+\alpha_1 Q'\big)
+\omega_2 \alpha_0 Q'\,.
\end{multline}
The first equation of (\ref{v2}) can be solved if the right hand side is is orthogonal 
to $Q'$ (the kernel of $L^+$). This imposes that either $\alpha_0=0$ or
$$
\int_{-\infty}^{\infty}\Big(-\alpha_0 Q'(x)
-\omega_1\Big(-\frac{\omega_1\alpha_0}{2}(xQ(x))
+\beta_1 Q(x)\Big)
-\omega_2 \beta_0 Q(x)\Big)Q'(x)dx=0,
$$
which implies that $\omega_1^2=-4\theta^{2}$, where 
$\theta\equiv\|Q'\|_{L^2(\mathbb{R})}/\|Q\|_{L^2(\mathbb{R})}$, i.e. $\omega_1=\pm i\theta$.
Hence if $\alpha_0\neq 0$ we have an eigenmode with purely imaginary $\omega_1$.

The second equation of (\ref{v2}) can be solved only if the right hand side is orthogonal to the
kernel of $L^-$, i.e. to $Q$. This imposes that either $\beta_0=0$ or
$$
\int_{-\infty}^{\infty}
\Big(-\beta_0 Q(x)
+\omega_1\big(\frac{\omega_1\beta_0}{2}(xQ'(x)+Q(x))+\alpha_1 Q'(x)\big)
+\omega_2 \alpha_0 Q'(x)\Big)
Q(x)dx=0
$$
which implies that $\omega_1^2=4$, i.e. $\omega_1=\pm 2$. 
From the above discussion, we have that either $\alpha_0=0$ or $\beta_0=0$.
If $\alpha_0\neq 0$ (and thus $\beta_0=0$) we obtain purely imaginary $\omega_1$ and have the bifurcation of
$(Q',0)$. These modes are not of interest for us.
If $\beta_0\neq 0$ (and thus $\alpha_0=0$) 
we indeed have en eigenvalue with positive $\omega_1$.
This mode corresponds to the eigenvector which is 
the bifurcation of $(u,v)=(0,Q)$ to the unstable mode of the form (\ref{form}) 
for the linearized about $Q$ cubic NLS equation.

%%%%%%%%%%%%%%%%%%%%%%%%%%%%%%%%%%%%%%%%%%%%%%%%%%%%%%%%%%%%%%%%%%%%%%%%%%%%%%%%%%%%%%%%%%%%%%%%%%%%%%%%%%
%%%%%%%%%%%%%%%%%%%%%%%%%%%%%%%%%%%%%%%%%%%%%%%%%%%%%%%%%%%%%%%%%%%%%%%%%%%%%%%%%%%%%%%%%%%%%%%%%%%%%%%%%%

\end{document}